\documentclass[12pt]{amsart}
\usepackage{amscd,amssymb,epsfig}
\usepackage{graphics}

\input epsf
\newtheorem{theorem}{Theorem}[section]
\newtheorem{lemma}[theorem]{Lemma}
\newtheorem{proposition}[theorem]{Proposition}
\newtheorem{corollary}[theorem]{Corollary}
\theoremstyle{definition}

\newtheorem{example}[theorem]{Example}

\newtheorem{question}[theorem]{Question}
\theoremstyle{remark}

\numberwithin{figure}{section}
\numberwithin{table}{section}

\def\M1g{{\mathcal M}_{g,1}}
\def\I1g{{\mathcal I}_{g,1}}

\catcode`\@=\active % make @ other character ( is default)

\begin{document}

\title[Stable curves and screens]
{Stable curves and screens on fatgraphs}
\author{R. C. Penner}
\address{Departments of Mathematics and Physics/Astronomy\\
University of Southern California\\
Los Angeles, CA 90089\\
USA\\
~{\rm and}~Department of Mathematics\\
Aarhus University\\
DK-8000 Aarhus C, Denmark\\}
\email{rpenner{\char'100}usc.edu}

\author{Greg McShane}
\address{Laboratoire Emile Picard\\
Universite Paris Paul Sabatier, UFR MIG \\
118 route de Narbonne\\
31062 Toulouse cedex 4\\
France}
\email{greg.mcshane{\char'100}gmail.com}

\keywords{moduli space of curves, stable curves, Deligne-Mumford compactification}
\subjclass{Primary 32G15, 57M99; Secondary 14H10, 14G15, 57N05, 20F99}

\thanks{RCP is happy to acknowledge useful discussions with Alex Bene, Kevin Costello, 
and Dennis Sullivan and to thank the Laboratoire Emile Picard in Toulouse and especially the Max Planck Institute for Mathematics 
in Bonn for warm hospitality.}

\begin{abstract}
The mapping class group invariant ideal cell decomposition of the Teichm\"uller space
of a punctured surface times an open simplex has been used in a number of computations.
This paper answers a question about the asymptotics of this decomposition, namely,
in a given cell of the decomposition, which curves can be short?  Screens are a new 
combinatorial structure which provide an answer to this question.  The heart of the calculation
here involves Ptolemy transformations and the triangle inequalities on lambda lengths.
\end{abstract}

\maketitle

\section{Introduction}

Throughout this paper, $F=F_g^s$ will denote
a fixed smooth oriented surface of genus $g$ with $s\geq 1$ punctures, where $2g-2+s>0$,
with mapping class group $MC(F)$.

Let ${\mathcal T}(F)$ denote the Teichm\"uller space of $F=F_g^s$ and
$\widetilde {\mathcal T}(F)$ denote the trivial ${\mathbb R}^s_{>0}$ bundle over it.
Let  ${\mathcal M}(F)={\mathcal T}(F)/(MC(F)$ denote Riemann's moduli space
with its Deligne-Mumford compactification $\overline{\mathcal M}(F)$.
$MC(F)$ also acts on $\widetilde {\mathcal T}(F)$ by permuting the numbers assigned to
punctures.

There is a $MC(F)$-invariant ideal cell decomposition \cite{Harer85,Hubbard-Masur,Penner87,Strebel} of $\widetilde{\mathcal T}(F)$ which has found wide application in geometry and physics \cite{AC,Harer86,HarerZagier,Igusa04,Kontsevich92,Mondello,MoritaPenner, Penner88}.  Cells in this
decomposition are in one-to-one correspondence with homotopy classes of ``fatgraph spines''  
of $F$, that is, a homotopy class of embedded graph in $F$ in the usual sense  together with
cyclic orderings on the half-edges about each vertex.  (See the next section 
for further precision)

Thus, to each fatgraph spine $G$ of $F$, there is a corresponding cell ${C}(G)\subset\widetilde{\mathcal T}(F)$.  In the interests of understanding $\overline{\mathcal M}(F)$ combinatorially, it is natural ask:  

\vskip .2in

\begin{question}\label{quest}
Given $G$ and given a collection $K$ of non-parallel and non-puncture-parallel disjointly embedded
and essential simple closed curves in $F$, when is there $\tilde\Gamma\in C(G)$ so that
the curves in $K$ are the only curves of hyperbolic length less than $\epsilon$ for the underlying hyperbolic metric, for some small $\epsilon >0$?  That is, which curves can be short in $C(G)$?
\end{question}

\vskip .2in

We give in this paper a complete answer to this question, as follows, where we shall concentrate
in this introduction on the case that $G$ is trivalent for simplicity.  

Let $E$ denote the set of 
edges of $G$ and consider any proper subset $A\subset E$.  There is a smallest 
(not necessarily connected) subgraph 
$G_A$ of $G$ containing $A$, and we say that $A$
 is ``recurrent'' if $G_A$ has no univalent vertices.
 (Again, see the next section
for a more detailed discussion of recurrence.)  
Suppose $A$ is recurrent and $G_A$ is connected,
and get rid of all bivalent vertices of
$G_A$ in the usual way to produce either
a simple cycle  in $G$ or another trivalent fatgraph $G'$.
A neighborhood of $G_A\subset G$ in $F$ is a subsurface of $F$, an annulus in the former case
and a punctured surface of negative Euler characteristic in the latter.  Define the ``relative boundary''
of $A$ to be the edge-path in $G$ of the simple cycle itself in the former case and those
of the boundary components of this subsurface in the latter case, where you discard any
such cycles that are puncture-parallel in $F$ itself.

The new combinatorial structure which provides the answer to Question~\ref{quest} (and was
introduced in \cite{Penner06a}), a 
``screen on a fatgraph $G$'' is a subset ${\mathcal A}$ of the power set (i.e., the set of subsets)
of the set $E$ 
of edges of $G$ with the following properties:

\vskip .2in

\leftskip .2in

\noindent {\rm i)}~$E\in{\mathcal A}$;

\vskip .1in

\noindent {\rm ii)}~each $A\in{\mathcal A}$ is recurrent;

\vskip .1in

\noindent {\rm iii)}~if $A,B\in{\mathcal A}$ with $A\cap B\neq\emptyset$, then either
$A\subseteq B$ or $B\subseteq A$;

\vskip .1in

\noindent {\rm iv)}~for each $A\in{\mathcal A}$, 
$\cup\{ B\in{\mathcal A}:B~{\rm is~a~proper~subset~of}~A\}$ is a proper subset of $A$. 

\vskip .2in

\leftskip=0ex

\noindent Condition i) is simply a convenient convention, conditions iii-iv) are familiar
from Fulton-MacPherson \cite{FultonMacPherson}, and here we impose the further condition ii) of recurrence. 
Notice that the properness condition iv) and recurrence condition ii) together imply 
that if $G_A$ is a simple cycle in $G$, then for any screen ${\mathcal A}$ on $G$ 
with $A,B\in{\mathcal A}$ and $A\cap B\neq\emptyset$, we must have $A\subseteq B$,
i.e., simple cycles are necessarily atomic in any screen.

Each element $A\in{\mathcal A}$ other than $A=E$ has an immediate predecessor
$A' \in {\mathcal A}$, and regarding $A$ as a set of edges
in $G_{A'}$ in the natural way, has its relative boundary $\partial _{\mathcal A} A$ defined before.
Finally, the ``boundary'' of the screen itself is 
$\partial {\mathcal A} = \bigcup_{A\in{\mathcal A}-\{ E\}}~\partial _{\mathcal A} A$.

Here is the answer to Question~\ref{quest}, our main result:

\vskip .2in

\begin{theorem} \label{mainthm} For any fatgraph $G$, the cell ${C}(G)$ admits as short curves a family $K$ of non-parallel and non-puncture parallel disjointly embedded
and essential simple closed curves in $F$ if and only if $K=\partial {\mathcal A}$ for some screen
${\mathcal A}$ on $G$.
\end{theorem}

\vskip .2in

Let us immediately do several examples, where it is typically easiest to study the quotient
of $\widetilde{\mathcal T}(F)$ by the natural ${\mathbb R}_{>0}$-action, the projectivized space, which we shall
denote
$$P\widetilde{\mathcal T}(F)=\widetilde{\mathcal T}(F)/{\mathbb R}_{>0}\approx {\mathcal T}(F)\times\Delta^{s-1},$$
where $\Delta ^p$ denotes the open $p$-dimensional simplex.
In particular for a once-puntured surface, we have
$P\widetilde{\mathcal T}(F)={\mathcal T}(F)$.

\vskip .2in

\centerline{\epsffile{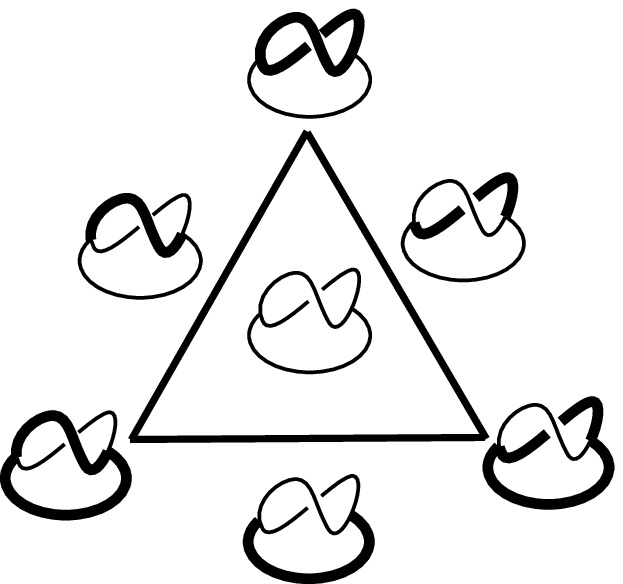}}
\vskip .1in

\centerline{{\bf Figure 1}~Screens on the once-punctured torus.}

\vskip .2in

\begin{example}
For the once-punctured torus $F=F_1^1$, the ideal cell decomposition of
${\mathcal T}(F)$ is the Farey tesselation of the disk \cite{Penner87}.  In Figure~1 we depict a typical
top-dimensional 2-cell, which is indexed by a non-planar fatgraph $G$ with two 3-valent vertices
as is also illustrated.
The codimension-one cells arise by collapsing any one of the three edges shown as darkened in the
figure, and
the codimension-two cells at infinity are indexed by the three possible recurrent
subgraphs of $G$ as likewise illustrated.  In this example, non-empty screens are always singletons, and the boundary of a screen always consists of a single curve.
\end{example}

\vskip .2in

\centerline{\epsffile{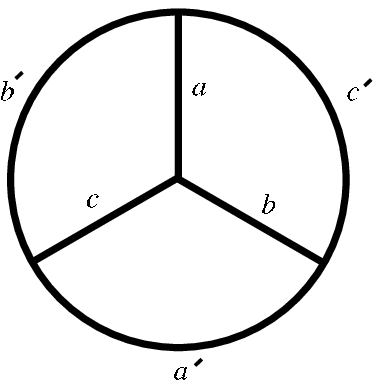}}
\vskip .1in

\centerline{{\bf Figure 2}~A fatgraph for the four-punctured sphere.}

\vskip .2in

\begin{example}
For the four-times punctured sphere $F=F_0^4$, consider the Mercedes sign fatgraph $G$ depicted
in Figure 2.  Both screens
${\mathcal A}_1=\{ \{ a,b,a',b'\}\}$  and
${\mathcal A}_2=\{\{ a,b,c,a',b'\} ,\{a,c,b'\}\}$
correspond to pinching to zero the closed edge-path $a-b-a'-b'$,
and both screens have this same edge-path as boundary.
\end{example}

\vskip .2in

\centerline{\epsffile{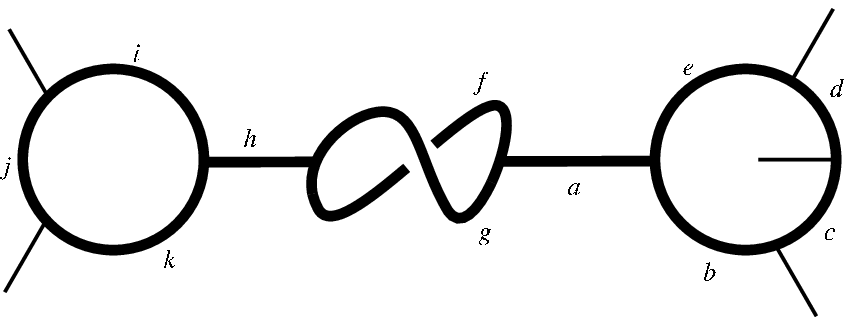}}
\vskip .1in

\centerline{{\bf Figure 3}~Typical example.}

\vskip .2in

\begin{example}
Consider the sub-fatgraph of a fatgraph $G$ with edges $E$ depicted in Figure~3 and the screen
$${\mathcal A}=\bigl\{ E,
\{ a-k\} , \{ a-g\} , \{ f-g\}  ,\{ b-e\} , \{ i-k\}
\bigr\}$$
on $G$.  The boundary of $\mathcal A$ is comprised of the four edge-paths
$f-g$, $b-c-d-e$, $h-i-j-k-h-f-g$, and $a-b-c-d-e-a-f-g$. 
\end{example}

\vskip .2in

We shall rely on ``lambda length''  coordinates from \cite{Penner87}  (recalled in $\S$3) on the {\it decorated Teichm\"uller space}
$\widetilde{\mathcal T}(F)$, where the fiber over a point
is taken to be the set of all $s$-tuples of horocycles, one horocycle about each puncture; one may take the hyperbolic lengths of the distinguished horocycles as a convenient coordinate on the fiber.  

In effect, we shall record the rates of divergence of lambda lengths regarded as projective
coordinates on $P\widetilde{\mathcal T}(F)$, and the
crucial point is that in the cell ${C}(G)$, the lambda length coordinates on $G$
must satisfy all three strict triangle inequalities at each vertex of $G$ (cf. Lemma \ref{triangle}).
This is what forces the recurrence condition.

The proof of Theorem~\ref{mainthm} involves the explicit calculation of holonomies
using  ``path-ordered products'' of matrices (due to Bill Thurston and Volodya Fock  \cite{fock} independently and recalled in $\S$3).   The proof further requires
estimates on the absolute traces of the representing matrices.  To this end, we find a
condition weaker than the triangle inequalities which satisfies two properties:
1) the condition is invariant under certain ``Whitehead moves'' (see the next section for a definition)
sufficient to simplify the path-ordered product; and 2) the condition guarantees
the required estimates on the absolute traces.  This is the heart of the paper (in $\S$5),
and the techniques involve only ``Ptolemy transformations'' 
(cf. Lemma~\ref{horos}a), path-ordered products, and the triangle inequality.

Because the argument at heart only depends upon these formulae, we are optimistic that the current paper may have ramifications more generally for cluster algebras \cite{GSV} and cluster ensembles
 \cite{FockGoncharov}.
Since Wolpert has recently announced \cite{Wolpert} that lambda lengths are strictly convex  along Weil-Peterssen
geodesics, we are likewise optimistic about applications to the asymptotic WP geometry.

There is furthermore a program to extend the cell decomposition of moduli space to the
Deligne-Mumford compactification using screens, which is already well underway
(as discussed in the closing remarks $\S$7).

\vskip .2in

\section{Fatgraphs and recurrence}

A {\it graph} is a finite one-dimensional CW complex with no isolated vertices
whose 1-cells are {\it edges} and whose 0-cells are {\it vertices}.  The set of edges of $G$
will be denoted $E(G)$.  A {\it half-edge}
of an edge $e\in E(G)$ is either one of the two components of the interior of $e$ with an interior
point removed, and the {\it valence} of a vertex is the number of half-edges containing the vertex
in their closures, said to be {\it incident} on the vertex.
A {\it fatgraph} is a graph together with a cyclic ordering on the half-edges incident on each
vertex.  In particular, a finite CW decomposition of a circle is an example of a fatgraph,
as is a planar tree where the cyclic ordering is induced by the counter-clockwise 
orientation on the plane.

A fatgraph $G$ determines
a punctured surface $F'(G)$ gotten by assigning to each $k$-valent vertex an oriented ideal
$k$-gon, whose sides correspond to the incident half-edges, and finally
identifying in the natural way pairs of sides of these polygons associated to pairs of half-edges contained in a common edge of $G$.  
 
The vertices of the ideal polygons are identified to the punctures of $F'(G)$.  Each edge $e\in E(G)$ gives rise to its dual ideal arc $\alpha _{(G,e)}$ connecting punctures 
in $F'(G)$.

An {\it ideal triangulation} of $F_g^s$ is the homotopy class of a set of arcs
connecting punctures in $F_g^s$, called {\it ideal arcs}, which decompose the surface into a collection of triangles with vertices at the punctures.  More generally, an {\it ideal cell decomposition} is the homotopy class of a subset of an ideal triangulation which decomposes the surface into polygons.  

Provided each vertex of $G$ has valence at least three, $\{ \alpha _{(G,e)}:e\in E(G)\}$ is an ideal cell decomposition of $F'(G)$ said to be {\it dual} to $G$.  Conversely,
the Poincar\'e dual of an ideal cell decomposition of $F_g^s$ is a fatgraph embedded in $F_g^s$ each of whose vertices has valence at least three, where the cyclic ordering in the fatgraph structure
is induced by the clockwise order in the oriented surface $F_g^s$.  
 
A fatgraph $G$ also determines a corresponding oriented surface $F(G)$ with boundary
constructed by assigning to each $k$-valent vertex an oriented $(2k)$-gon, whose alternating
sides correspond to the incident half-edges, and as before, identifying pairs of sides of these
polygons corresponding to pairs of half-edges contained in a common edge of $G$.  
The alternating unpaired edges of these polygons comprise the boundary of $F(G)$.  We may 
regard $F(G)\subseteq F'(G)$ as a strong deformation retraction in the natural way.

In particular, $G$ is a strong deformation retraction or {\it spine} of $F(G)$ or $F'(G)$.  It follows that any free homotopy class
of essential curve in $F(G)$ or $F'(G)$ gives rise to a closed edge-path in $G$, which is uniquely
determined up to its starting point provided we demand that the edge-path is {\it efficient}
in the sense that it never consecutively traverses the same edge of $G$ with opposite
orientations.  

A closed edge-path in $G$ corresponding to a boundary component of $F(G)$
will be called a {\it boundary component} of $G$ itself.  An efficient 
boundary component of $G$ must have edges of $G$ incident on only one side.
Put another way for a trivalent fatgraph, an efficient edge-path is a boundary component
if and only if it consists entirely of left turns or consists entirely of right turns.

Suppose that $G$ is a fatgraph with set $E=E(G)$ of edges and corresponding surface $F=F(G)$.  Any subset $A\subseteq E$
determines a subgraph by including all vertices of $G$ on which edges in $A$ are incident.
Furthermore by restriction, the fattening on $G$ induces a fattening on this subgraph, which
thus determines a well-defined sub-fatgraph $G_A$.  We may regard $F(G_A)$ as a subsurface
embedded in the interior of $F$ in the natural way.
Define the {\it boundary} of $A$ to be the collection $\partial A$
of (unoriented) efficient closed edge-paths in $G$ corresponding to 
the relative boundary of $F(G_A)$ in $F=F(G)$, that is,
the collection of closed edge-paths corresponding to the components of 
the boundary $\partial F(G_A)$ which are not homotopic to boundary components
of $F$ itself.  In particular, if $G_A$ is a circle, then $\partial A$ 
is the closed edge path of $G_A$ if this circle is not boundary parallel in $F$,
and $\partial A$ is empty if this circle is boundary parallel in $F$.

We say that $A\subseteq E$ is {\it recurrent} if for every edge $a\in A$, there is
an efficient closed edge-path $\gamma _a$ in $G$ so that $\gamma _a$ traverses $a$
and traverses only edges in $A$.  Any subset $A\subseteq E$ has a (possibly empty, e.g., in
the case of a planar tree)
{\it maximal recurrent subset} $R(A)$, namely, the set of edges of $A$ traversed by
an efficient closed edge-path in $G_A$.

\begin{lemma}\label{recurrence}
Suppose that $G$ is a fatgraph and $A\subseteq E=E(G)$.
Then the following are equivalent:

\vskip .2in

\leftskip .2in

\noindent{\rm i)} $A$ is recurrent;

\vskip .1in

\noindent {\rm ii)}
there is a function $\mu: E\to {\mathbb Z}_{\geq 0}$  whose support is $A$ so that
for each vertex of $G$ with incident half-edges $e_1,\ldots ,e_k$ and extending
the function $\mu$ to be defined on half-edges in the natural way,
we have that  $\sum _{i=1}^k \mu (e_i)$ is even,
and the generalized weak triangle inequalities hold, i.e., for
each $j=1,\ldots ,k$,
$$\mu (e_j)\leq \sum _{i\neq j}
\mu (e_i);$$

\vskip .1in

\noindent {\rm iii)} every vertex of $G_A$ has valence at least two.

\leftskip=0ex

\end{lemma}

\vskip .2in

\begin{proof} First suppose that $A$ is recurrent, and let
$\mu _a(e)$ be the number of times that a chosen $\gamma_a$ traverses $e$ for each
$a\in A$ and $e\in E$.  
Each $\mu _a: E\to {\mathbb Z}_{\geq 0}$ satisfies the restrictions of condition (ii), hence so too does their sum $\mu =\sum _{a\in A} \mu _a$, which
has full support on $A$.  Thus, (i) implies (ii).

Conversely, suppose that  $\mu$ is a function supported on $A$ satisfying the properties
of condition (ii).  For each $k$-valent vertex of $G$, there is a dual ideal
$k$-gon
in the corresponding punctured surface $F'(G)$, and we shall construct a family
of arcs properly embedded in this $k$-gon realizing the values of $\mu$ on
the dual edges of $G$ as the geometric intersection numbers.
These arc families in the $k$-gons then combine uniquely to produce disjointly embedded
curves in the natural way, 
whose component simple closed curves 
in $F$ have corresponding edge-paths  which satisfy the required properties.

The construction in each $k$-gon
proceeds by induction on $k\geq 2$ with notation for incident edges as in condition (ii).  In case
$k=2$, simply take a collection of $\mu(e_1)=\mu(e_2)$ arcs crossing the bigon.  For the 
case $k=3$, take $${1\over 2} [\mu (e_{i_1})+\mu (e_{i_2})-\mu (e_{i_3})]
={1\over 2} [\mu (e_{i_1})+\mu (e_{i_2})+\mu (e_{i_3})-2\mu (e_{i_3})]~\in ~{\mathbb Z}_{\geq 0}$$ parallel
copies of the arc joining edges $e_{i_1}$ to $e_{i_2}$, where $\{ i_1,i_2,i_3\}=\{ 1,2,3\}$.  For the induction step, take a consecutive pair of edges
$e_i,e_{i+1}$ so that $\mu (e_i)+\mu (e_{i+1})$ is least among all consecutive
pairs of edges, here taking the indices modulo $n$ so that $e_{n+1}=e_1$.
Cutting along the diagonal separating $e_i$ and $e_{i+1}$ from the rest 
decomposes the $k$-gon into a $(k-1)$-gon and a triangle. 
Extend $\mu$ to a function defined on the edges of these regions by taking
value $\mu (e_i)+\mu (e_{i+1})$ on the diagonal, so the generalized triangle
inequalities hold on each region by our choice of consecutive edges, and the
parity condition holds by construction.  By
the inductive hypothesis, appropriate arc families exist in each region, and
they combine in the natural way to give the required arc family in the $k$-gon
itself.  It follows that (i) is equivalent to (ii).

If $G_A$ has a univalent vertex, say with incident edge $a\in A$, then there can
be no efficient edge-path in $G_A$ traversing $a$, so (i) implies (iii).  To see that
(iii) implies (ii), define $\mu$ to take value 2 on the edges in $A$ and vanish otherwise,
and note that $\mu$ satisfies condition (ii) provided $G_A$ has no univalent
vertices.
\end{proof}
\vskip .2in

\vskip .2in

Suppose that $e$ is an edge of a fatgraph $G$ with distinct endpoints.  We may {\it collapse} $e$
to a vertex to produce a new fatgraph  $G'$, where the cyclic ordering at the resulting vertex
arises by combining the cyclic orderings on the half-edges incident on $e$ in the natural way.
Dually, one removes the ideal arc $\alpha _{(G,e)}$ from the dual ideal cell decomposition.

If $G$ is a trivalent fatgraph and $e$ is an edge of $G$ with distinct endpoints, then
a {\it Whitehead move}  on $e$ is the fatgraph that results by collapsing $e$ and then
un-collapsing the resulting four-valent vertex in the unique distinct manner.  A Whitehead move
along an edge $e$ is depicted in Figure~4, which furthermore indicates the notation near an
edge $e$ which we shall adopt in many of the calculations of this paper.

\vskip .2in

\centerline{\epsffile{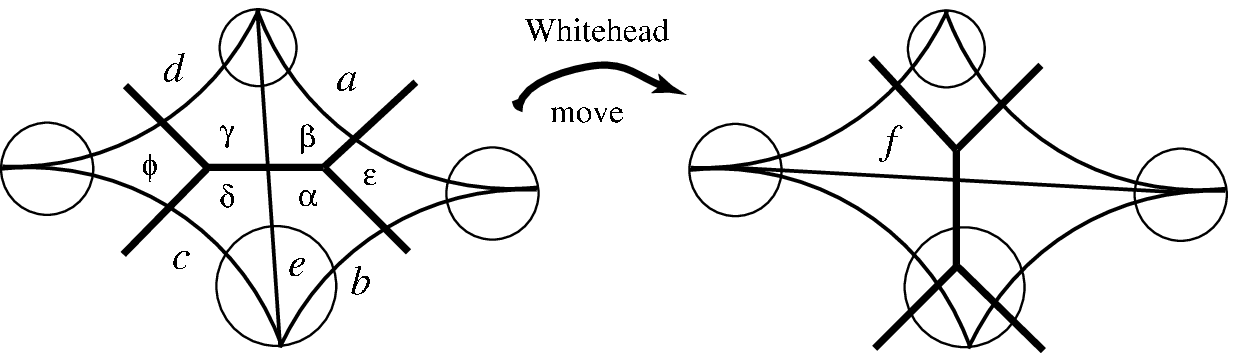}}
\vskip .1in

\centerline{{\bf Figure 4}~Standard notation for Whitehead moves.}

\vskip .2in

Using the characterization Lemma~\ref{recurrence}iii), it follows directly that
recurrence is invariant under Whitehead moves on trivalent fatgraphs and is furthermore in any case invariant under collapse of edges $e$ with distinct endpoints neither of which is univalent.

\vskip .2in

\section{Coordinates}

The reader is referred to \cite{Penner87} or the more recent treatment \cite{Penner06} for proofs and
further details on the material which is recalled in this section.  
We begin with several formulae on horocycles in the hyperbolic plane.

 If $h,h'$ are horocycles in the hyperbolic plane with distinct centers in the circle at infinity, then consider the unique geodesic $\gamma (h,h')$ connecting their centers.  The horocycles $h,h'$ truncate $\gamma(h,h')$ to a geodesic segment of
some finite signed length taken to be positive if and only if $h$ and $h'$ are disjoint.  Define
the {\it lambda length} of $h,h'$ to be $\lambda (h,h')=\sqrt{{\rm exp}~\delta}$.  (This is a different
normalization for lambda lengths than in \cite{Penner87}, for instance, where the lambda length
is taken as $\sqrt{2~{\rm exp}\delta}$, cf. \cite{Penner06}.)

\vskip .2in

\begin{lemma}\label{horos}
Suppose $h_1,h_2,h_3,h_4$ are horocycles with distinct centers occurring in this
clockwise order in the circle at infinity, and let $\lambda_{ij}=\lambda(h_i,h_j)$ for distinct
$i,j\in\{ 1,2,3,4\}$.  Then:

\vskip .2in

\leftskip .2in

\noindent {\bf a)}~{\rm [Ptolemy's equation]}~$\lambda_{13}\lambda_{24}=\lambda_{12}\lambda_{34}+\lambda_{14}\lambda_{23}$;

\vskip .1in

\noindent{\bf b)}~{\rm [Cross Ratios]}~the M\"obius transformation that takes the centers of $h_3,h_2,h_1$ respectively to $0,1,\infty$ also maps the center of $h_4$ to $-{{\lambda_{23}\lambda_{14}}\over{\lambda_{12}\lambda_{34}}}$;

\vskip .1in

\noindent{\bf c)}~{\rm [h-lengths]}~the hyperbolic length of the horocyclic segment in $h_i$ with endpoints
$h_i\cap \gamma(h_i,h_j)$ and $h_i\cap\gamma(h_i,h_k)$ is given by 
${{\lambda_{jk}}\over{\lambda_{ij}\lambda_{ik}}}$, where $\{ i,j,k\} =\{ 1,2,3\}$;

\vskip .1in

\noindent{\bf d)}~{\rm [Affine duality]}~taking the upper sheet ${\mathbb H}$ of the hyperboloid in Minkowski 3-space as the model for the hyperbolic plane, there is a unique isotropic vector $u_i$ with positive $z$-coordinate so that $h_i=\{ w\in{\mathbb H}: w\cdot u_i=-2^{-{1\over 2}}\}$, for $i=1,2,3,4$ where $\cdot$ denotes the pairing with quadratic form
$x^2+y^2-z^2$, and $\lambda_{ij}=\sqrt{-u_i\cdot u_j}$ for distinct $i,j\in\{ 1,2,3,4\}$;

\vskip .1in

\noindent{\bf e)}~{\rm [Simplicial coordinates]} in the notation of part d), the signed volume of the Euclidean tetrahedron in Minkowski 3-space spanned by
$u_1,u_2,u_3,u_4$ is given by $2\sqrt{2}~\lambda_{12}\lambda_{23}\lambda_{34}\lambda_{14}$ times
$${{\lambda_{12}^2+\lambda_{23}^2-\lambda_{13}^2}\over{\lambda_{12}\lambda_{23}\lambda_{13}}}
+
{{\lambda_{14}^2+\lambda_{34}^2-\lambda_{13}^2}\over{\lambda_{14}\lambda_{34}\lambda_{13}}},$$
where the sign is positive if and only if the edge of the tetrahedron connecting $u_1,u_3$ lies below the
edge connecting $u_2,u_4$.

\vskip .1in

\noindent{\bf f)}~{\rm [Ellipticity]}~in the notation of part d), the affine plane containing $u_1,u_2,u_3$
determines an elliptic conic section if and only if $\lambda_{12},\lambda_{13},\lambda_{23}$ satisfy
the three strict triangle inequalities.

\leftskip=0ex
\end{lemma}

Given a point $\tilde\Gamma\in\widetilde{\mathcal T}(F)$ and given the homotopy class of an ideal arc $\alpha$
in $F$, we may straighten $\alpha$ to the geodesic for the underlying
hyperbolic structure and truncate this geodesic by cutting it at the horocycles centered at its
endpoints coming from the decoration.  This geodesic segment has a signed hyperbolic length
$\delta$ taken with a positive sign if and only if the horocycles are disjoint.  The basic coordinate
of an ideal arc in a decorated hyperbolic surface is the {\it lambda length} (also sometimes
called the ``Penner coordinate'') defined by
$\lambda (\alpha;\tilde\Gamma)=\sqrt{{\rm exp}~\delta}$.  

\vskip .1in

\begin{theorem} \label{lambda}
Fix any trivalent fatgraph $G$.
Then the assignment of lambda lengths
$$\aligned
\widetilde{\mathcal T}(F'(G))&\to {\mathbb R}_{>0}^{E(G)}\\
\tilde\Gamma &\mapsto (e\mapsto \lambda (\alpha _{(G,e)};\tilde\Gamma))\\
\endaligned$$
is a real-analytic homeomorphism onto.
\end{theorem}

\vskip .2in

For convenience when the fatgraph $G$ is fixed or understood, we shall refer to
the lambda length of an edge $e$ of $G$ rather than that of its dual arc $\alpha _{(G,e)}$.
We shall also often identify an arc with its
lambda length for convenience.

Suppose that $G$ is a trivalent fatgraph.  Consider an edge $e$ of $G$ and adopt the notation of Figure~4, where $e$ has distinct endpoints with incident half-edges $a,b$ and $c,d$ occurring in the alphabetic counter-clockwise order about $e$.  (If $e$ does not have distinct endpoints or if $a,b,c,d$ are not distinct, then
adopt the corresponding notation for nearby edges in the universal cover.)  
Dual to each vertex of $e$ is an ideal triangle, and each such triangle has three vertices, denoted by
Greek letters in Figure~4.  To each such triangle/vertex pair is naturally associated a {\it sector} of
$G$, that is, a pair of consecutive half-edges of $G$ incident on a common vertex, namely, the pair of
half-edges adjacent to the given vertex in the given triangle.

In fact, one can conveniently calculate the holonomies of based closed curves in $F'(G)$ as follows.
Define the matrices
$$R=\left (\begin{array} {cc} \hskip 1.5ex 1&1\\ -1&0\\\end{array}\right ),~~
L=\left (\begin{array} {cc}0&-1\\ 1&\hskip 1.5ex 1\\\end{array}\right )~~\in PSL_2({\mathbb R}).$$
 According to Lemma~\ref{horos}b),
the cross ratio of the ideal quadrilateral with edges $\alpha _{(G,a)}$, $\alpha _{(G,b)}$, $\alpha _{(G,c)}$, $\alpha _{(G,d)}$ is given by $-bd/ac$, where we have identified an edge of $G$ with its lambda length for convenience, and we further define the matrix
$$X_e=\left (\begin{array} {cc}0&\sqrt{{{ac}/{bd}}}\\ -\sqrt{{bd}/{ac}}&0\\\end{array}\right ).$$

Choosing a vertex of $G$ as basepoint, consider a closed edge-path $\gamma $ in $G$ representing an essential based closed curve in $F$.  We may as well assume that $\gamma$
is efficient (though this is not necessary since $RL=R^3=L^3=X_e^2=1\in PSL_2({\mathbb R})$), so that it alternately traverses edges and sectors of $G$ and makes
turns, right or left, at each sector.  Suppose that $\gamma$ serially makes turns $t_i$ at the sectors, then traverses edges $e_i$, for $i=1,\ldots ,n$, and associate the {\it path-ordered product}
$$M=T_1X_{e_1}T_2\cdots X_{e_n}$$ of matrices,
where $T_i=R$ or $L$ if $t_i$ is a right or left turn respectively.  The matrix $M\in PSL_2({\mathbb R})$
gives the holonomy of the based curve $\gamma$.  Of course by conjugacy invariance of trace,
the absolute value of the trace of $M$ is independent of the basepoint.  We shall use these path-ordered products to detect the short curves that occur on a path in ${C}(G)\subseteq \widetilde {\mathcal T}(F'(G))$.

The quadrilateral in Figure~4 is realized as a geodesic ideal quadrilateral with horocycles centered at each vertex.  We define the {\it h-length} of a sector of $G$
to be the hyperbolic length of the corresponding horocyclic segment.
According to Lemma~\ref{horos}c), the h-length of a sector is the opposite lambda length divided by the product of adjacent lambda lengths.

Furthermore in the notation of Figure~4, we define the {\it simplicial coordinate} of the edge $e$ to be the
quantity
$${{a^2+b^2-e^2}\over{abe}}+{{c^2+d^2-e^2}\over{cde}}={a\over{be}}+{b\over{ae}}-{e\over{ab}}
+{c\over{de}}+{d\over{ce}}-{e\over{cd}}.$$
According to Lemma~\ref{horos}e), the simplicial coordinate is a multiple of the
signed volume of the corresponding tetrahedron, and by inspection, it is a linear combination
of the nearby h-lengths.  From the definition, the simplicial coordinate is the sum of two terms
each of which is associated to a vertex of the corresponding edge.

Consider a trivalent fatgraph $G$ with set $E$ of edges for the surface $F_g^s$ together with an assignment of lambda lengths $\lambda :E\to{\mathbb R}_{>0}$.  We say that $\lambda$ satisfies
the {\it no vanishing cycle condition} provided that all the corresponding simplicial coordinates are non-negative and there is no cycle in $G$ all of whose simplicial coordinates vanish.

\vskip .2in

\begin{theorem}\label{decomp}  For any surface $F=F_g^s$ with $s\geq 1$,
there is a $MC(F)$-invariant ideal cell decomposition of $\widetilde{\mathcal T}(F)$, where the
cells in this decomposition are in one-to-one correspondence with homotopy classes of embeddings
of fatgraph spines of $F$ each of whose vertices has valence at least three.  The face relation in this cell decomposition is generated by
Whitehead collapse.

In particular, if $G$ is a trivalent fatgraph spine
of $F$, then the closed cell $C(G)\subseteq \widetilde{\mathcal T}(F)$ corresponding to it
is described in lambda length coordinates with respect to $G$ by the no vanishing cycle condition.

Furthermore, suppose $G'$ arises from $G$ by collapsing to a point each component of a forest in $G$. Then the corresponding closed cell
$C(G')\subseteq C(G)\subseteq \widetilde{\mathcal T}(F)$ is described by taking all simplicial
coordinates on edges in the forest to vanish.  
Finally, any assignment of non-negative real numbers to the
edges of $G'$ with no vanishing cycles
is realized as the simplicial coordinates of a uniquely determined collection of positive lambda
lengths on $G$.

\end{theorem}~

~~\vskip .2in

\begin{corollary}
There is a $MC(F)$-invariant ideal cell decomposition of $P\widetilde{\mathcal T}(F)$ for any surface
$F=F_g^s$ with $s\geq 1$, where the
cells in this decomposition are in one-to-one correspondence with homotopy classes of embeddings
of fatgraph spines of $F$ whose vertices have valence at least three.
\end{corollary}~

\vskip .2in

\begin{proof}This follows immediately from the previous theorem and homogeneity
of the formula for simplicial coordinates.
\end{proof}

\vskip .2in

\begin{lemma}~~\label{tele}
Suppose that $\gamma$ is an efficient edge-path in $G$
serially traversing edges $e_i$ alternating with sectors $t_i$, for $i=1,\ldots , n$.  Let $E_i$ denote the
simplicial coordinate of $e_i$ and $\alpha _i$ the h-length of the sector $t_i$.  Then
$\sum _{i=1}^n E_i = 2\sum _{i=1}^n \alpha _i$.
\end{lemma}

\vskip .2in

\begin{proof}
The proof follows from the definition of simplicial
coordinates in terms of h-lengths.
\end{proof}~

~~\vskip .2in

\begin{lemma}~{\cite{Penner87}}~\label{triangle}
The no vanishing cycle condition 
implies that the lambda lengths
at any vertex of $G$
satisfy the three strict triangle inequalities.
\end{lemma}

\vskip .2in

\begin{proof}
Adopt the notation of Figure 4  for the half-edges near an edge $e$ (again, in the universal
cover if the edges $a,b,c,d$ are not distinct or if $e$ does not have distinct endpoints).  If
$c+d\leq e$, then $c^2+d^2-e^2\leq -2cd$, so the non-negativity of the simplicial coordinate $E$ of $e$ gives
$0\leq cd  [(a-b)^2-e^2]$, and we find a second vertex so that the triangle inequality fails.
This is a basic algebraic fact about simplicial coordinates.
It follows that
if there is any such
vertex so that the triangle inequalities do fail for the lambda lengths of incident half-edges,
then there must be an efficient closed edge-path $\gamma$ passing through such triangles.
Letting $e_i$ denote the consecutive edges of $G$ serially traversed by $\gamma$ and 
$b_i$ denote the half-edge of $G$ incident on the common endpoint of $e_i$ and $e_{i+1}$, we find
$e_{j+1}\geq b_j+e_j$, for $j=1,\ldots n$.  Upon summing and canceling
like terms, we find $0\geq \sum _{j=1}^n b_j$, which is 
absurd since lambda lengths are positive.
\end{proof}

\vskip .2in

\section{Screens}

Suppose that $G$ is a trivalent fatgraph with set $E$ of edges and corresponding surface
$F$, and suppose that $\lambda _t:E\to {\mathbb R}_{>0}$, i.e., $\lambda _t\in {\mathbb R}_{>0}^E$,
is a continuous one-parameter family of lambda lengths for $t\geq 0$.  We shall typically 
apply Theorem \ref{lambda} to regard such a one-parameter family as a path in
$\widetilde{\mathcal T}(F'(G))$ itself.
There is an induced 
$\bar\lambda_t\in P({\mathbb R}_{>0}^E)$, where $P$ denote projectivization, and by compactness,
there is an accumulation point of ${\rm lim}_{t\to\infty} \bar\lambda _t$ in $P({\mathbb R}_{\geq 0}^E)$.

Say that $\lambda _t$ {\it stable} if there is a unique such limit point
denoted $\bar\lambda _\infty\in P({\mathbb R}_{\geq 0}^E)$.  If $\lambda _t$ is any path, then any accumulation point of $\bar\lambda_t$ is also the limit of some stable path since decorated
Teichm\"uller space is path connected.

Suppose that $\lambda _t\in P({\mathbb R}_{>0}^E)$ is stable with limit 
$\bar\lambda_\infty\in P({\mathbb R}_{\geq 0}^E)$.  Set $E^0=E$ and make the following
recursive definition for $k\geq 1$.  If $E^{k-1}\neq\emptyset$, then there is $e\in E^{k-1}$ so that
for all $f\in E^{k-1}$, we have $\bar\lambda_\infty (e)/\bar\lambda_\infty (f)<\infty$, and we
set 
$$E^k=\{ f\in E^{k-1}: \bar\lambda_\infty (e)/\bar\lambda_\infty (f)=0\} ,$$
and set $E_k=\emptyset$ if $E_k=E_{k-1}$.
Notice that if two edges $e_1,e_2\in E^{k-1}$ satisfy the condition that
for all $f\in E^{k-1}$, we have $\bar\lambda_\infty (e_i)/\bar\lambda_\infty (f)<\infty$, for $i=1,2$,
then we have $0<\bar\lambda _\infty (e_1)/\bar\lambda_\infty(e_2)<\infty$, so $e_1$ and $e_2$ determine the same subset $E^k$.
Thus, $E=E^0\supseteq E^1\supseteq\cdots\supseteq E^N$ is a well-defined nested
sequence of finite length $N$ of proper subsets, and we set $E^{N+1}=\emptyset$ for convenience.

Now, suppose that $\lambda _t$ stays for all finite
$t\geq 0$ in the closed cell $C(G)$ corresponding to $G$, i.e., the lambda lengths satisfy the no
vanishing cycle condition by Theorem~\ref{decomp}.  
Define $${\mathcal A}(\lambda _t)=\{ A\subseteq E:A~{\rm is~the~set~of~edges~of~a~component~of~some}~
{E^k}\} ,$$
a subset of the power set of $E$.

\vskip .2in

\begin{proposition}\label{screen}
For any connected trivalent fatgraph $G$ with set $E$ of edges and any continuous stable one-parameter family $\lambda _t\in{\mathbb R}_{>0}^E$, for $t\geq 0$, which stays for all finite 
$t$ in the cell $C(G)\subseteq \widetilde{\mathcal T}(F(G))$ corresponding to $G$, the collection ${\mathcal A}={\mathcal A}(\lambda _t)$ satisfies the following properties:

\vskip .2in

\leftskip .2in

\noindent {\rm i)}~$E\in{\mathcal A}$;

\vskip .1in

\noindent {\rm ii)}~each $A\in{\mathcal A}$ is recurrent;

\vskip .1in

\noindent {\rm iii)}~if $A,B\in{\mathcal A}$ with $A\cap B\neq\emptyset$, then either
$A\subseteq B$ or $B\subseteq A$;

\vskip .1in

\noindent {\rm iv)}~for each $A\in{\mathcal A}$, 
$\cup\{ B\in{\mathcal A}:B~{\rm is~a~proper~subset~of}~A\}$ is a proper subset of $A$. 

\vskip .1in

\leftskip=0ex
\end{proposition}

\vskip .2in

\noindent A subset of the power set of $E$ satisfying properties i-iv) is called a {\it screen}
on $G$ for any (not necessarily trivalent) recurrent fatgraph $G$ with set $E$ of edges.

\vskip .2in

\begin{proof}  The first condition holds since $G$ is connected and $E=E^0$.
Recursively applying Lemmas~\ref{recurrence} and \ref{triangle},
we conclude that $E^k$ is a proper recurrent set in the possibly disconnected
fatgraph $G_{E^{k-1}}$, for $k=1,\ldots ,N$, so the second condition holds as well.
The third condition holds since two components
of a topological space either coincide or are disjoint, and the fourth follows since each inclusion
$E^{k}\subseteq E^{k-1}$ is proper.
\end{proof}

\vskip .2in

If ${\mathcal A}$ is a screen, then each $A\in{\mathcal A}-\{ E\}$ has an immediate
predecessor $A'$, i.e., $A\subseteq A'$ and if $B\in{\mathcal E}$ and $A\subseteq B\subseteq A'$, then$B=A$ or $B=A'$.  The maximum length of a chain $A\subseteq A'\subseteq \cdots \subseteq E$
of immediate predecessors in ${\mathcal A}$ is called the {\it depth} of $A$ in ${\mathcal A}$,
and the depth of $e\in E$ in ${\mathcal A}$ is the maximum depth of $A\in{\mathcal A}$ with
$e\in A$.

\vskip .2in

\begin{lemma} \label{cubicrealization}
Every screen ${\mathcal A}$ on every trivalent fatgraph $G$ arises
as ${\mathcal A}={\mathcal A}(\bar\lambda _t)$ for some stable $\lambda _t\in{\mathbb R}_{>0}^E$
lying in $C(G)$.
\end{lemma}

\vskip .2in

\begin{proof}
For any screen ${\mathcal A}$ on any trivalent fatgraph $G$, define a one-parameter family of
lambda lengths by taking $\lambda _t(e)=t^{d_e}$, where $d_e$ is the depth of $e$ in
${\mathcal A}$.  For any vertex $v$ of $G$, the maximum degree of the incident (half-)edges 
is achieved either twice or thrice by recurrence of elements of ${\mathcal A}$.  Thus,
the contribution from $v$ to each of the three possible simplicial coordinates of edges incident on $v$
is positive, and so the simplicial coordinate of each edge of $G$ for $\lambda _t$ is also positive;
$\lambda _t$ thus lies in $C(G)$ by Theorem~\ref{decomp}, and ${\mathcal A}(\bar\lambda _t)={\mathcal A}$ by construction.
\end{proof}

\vskip .2in

Let $\partial _{\mathcal A} A$ denote the relative boundary of $F(G_A)$ in $F(G_{A'})$, where $A'$ is the immediate predecessor of $A$ in ${\mathcal A}$, and define the {\it boundary} of ${\mathcal A}$ itself
to be
$$\partial {\mathcal A} =\bigcup _{A\in{\mathcal A}-\{ E\}} \partial _{\mathcal A} A.$$
~

\vskip .2in

\begin{lemma}~\label{shortcurves}
For any trivalent fatgraph $G$ with set $E$ of edges and stable $\lambda _t\in{\mathbb R}_{>0}^E$ lying  in $C(G)$, each edge-path in $\partial {\mathcal A}(\bar\lambda _t)$ is homotopic to a curve in $F'(G)$
whose hyperbolic length tends to zero as $t$ tends to infinity.  Furthermore, these are only such asymptotically short curves
for $\lambda _t$.
\end{lemma}

\vskip .2in

\begin{proof}
Let $K$ be a component of $\partial {\mathcal A}$, so $K=\partial _{\mathcal A} A$ for some
$A\in{\mathcal A}-\{ E\}$ with immediate predecessor $A'$.
Orient $K$ with the subsurface $F(G_A)$ on its left.  Consider the
universal cover $\tilde F$ of $F=F(G)$, let $\tilde G$, $\tilde G_A$, and $\tilde G_{A'}$ respectively
denote the full pre-images of $G$, $G_A$, and $G_{A'}$ in $\tilde F$, and choose
a lift $\tilde K$ of $K$ to $\tilde F$.  We shall refer to lambda lengths of edges of $\tilde G$, by which we
mean the value of $\lambda_t$ on the projection of the edge to $F$, and we will as usual
identify the lambda length of an edge of $\tilde G$ with the edge itself for convenience.

On the right of $\tilde K$ since $K$ is homotopic to a 
boundary component of $F(G_A)$, there are no edges of $\tilde G_A$,
and since $K$ is not homotopic to a 
boundary component of $F(G_{A'})$, there is at least one edge of 
$\tilde G_{A'}$ not in $\tilde G_A$ on the right.
Furthermore on the left of $\tilde K$, there is at least one
edge of $\tilde G_{A'}$ again since $K$ is not homotopic to a boundary
component of $F(G_{A'})$.

Since 
$\lambda_t$ corresponds to points in $C(G)$, it follows that the triangle inequalities hold on lambda lengths at each vertex of $\tilde G$ by Lemma~\ref{triangle}.  We claim that the following further properties of lambda lengths follow from these facts, where all limits are taken as $t\to\infty$:
\vskip .2in

\leftskip .2in

\noindent 1) if $x$ is an edge on the right of $\tilde K$ and $y$ is an edge of $\tilde K$, then 
we have ${x\over y}\to 0$;

\vskip .1in

\noindent 2) if $x$ is an edge on the right of $\tilde K$, $y_0$ is an edge on the left of $\tilde K$, and $y_1,y_2$ are edges of $\tilde K$ so that $y_0,y_1,y_2$ are all incident at a common vertex in $\tilde K$, then
${{xy_0}\over{y_1y_2}}\to 0$;

\vskip .1in

\noindent 3) if $y_0,y_1$ are consecutive edges of $\tilde K$ and $x$ is an edge on the right
of $\tilde K$ incident on their common endpoint, then ${{y_0}\over{y_1}}\to 1$.

\vskip .2in

\leftskip=0ex

The first property follows from the definition of $K$ as a relative boundary component of
$F(G_A)$ in $F(G_{A'})$ and the definition of the screen ${\mathcal A}(\lambda _t)$.  For property 2,  $y_0,y_1,y_2$ satisfy the triangle inequality $y_0<y_1+y_2$, so dividing by $y_1y_2$ and multiplying by $x$, we find
${{xy_0}\over{y_1y_2}}< {x\over{y_1}}+{x\over{y_2}}$; the right hand side tends to zero
by property 1.  Finally for property 3 again by the triangle inequalities, we have
$y_0<y_1+x$ and $y_1<y_0+x$.  Upon dividing the first by $y_1$ and the second by $y_0$
and applying property 1, we conclude $1\leq {\rm lim} ~{{y_0}\over{y_1}} \leq 1$, as required.

The first key point about properties 1-3) is that they are invariant under certain Whitehead moves.
In each case, we shall perform a Whitehead move along an edge $e\in K$,
where one vertex of $e$ has incident half-edges $a,b$ and the other vertex has
incident half-edges $c,d$, and where the edges $a,b,c,d$ occur in this counter-clockwise order
about $e$.  We shall refer to properties 1-3) for the fatgraph before the Whitehead move
and the corresponding properties 1$'$-3$'$) for the resulting fatgraph, and we shall let $f={{ac+bd}\over e}$
denote the edge and lambda length of the edge resulting from $e$ under the Whitehead move.

The first case of utility is when $b,c,e$ lie in $\tilde K$ and $a,d$ lie on the right of $\tilde K$.
The properties for this fatgraph imply that: 1) ${x\over y}\to 0$ for $x\in\{a,d\}$ and $y\in\tilde K$;
2) does not involve the vertices of $e$; and 3) ${b\over e}\to 1$ and ${e\over c}\to 1$.
Property 3$'$) requires ${b\over c}\to 1$, which follows from property 3).
Furthermore by the Ptolemy equation,
$$\aligned
{f\over a} ={{ac+bd}\over {ae}}&={c\over e} + {b\over e}{d\over a}\to 1+{d\over a},\\
{f\over d} ={{ac+bd}\over {de}}&={b\over e} + {c\over e}{a\over d}\to 1+{a\over d},\\
\endaligned$$
since ${c\over e}\to 1$ and ${b\over e}\to 1$.  Thus, at least one of ${f\over a},{f\over d}$
has a finite limit, hence ${f\over y}={f\over a} {a\over y}={f\over d}{d\over y}\to 0$
for any $y\in\tilde K$ by property 1) proving property 1$'$) and likewise for property 2$'$).

The second case of utility is when $b,d,e$ lie in $\tilde K$ with $a$ on the right and
$c$ on the left of $\tilde K$.  The properties for this fatgraph imply that: 1) ${a\over y}\to 0$
for  any $y$ in $\tilde K$; 2) ${{xc}\over{de}}\to 0$ for any $x$ on the right; and
3) ${b\over e}\to 1$.  Property 3$'$), namely, ${d\over f}\to 1$, follows from
$${f\over d}={{ac+bd}\over{de}}={{ac}\over{de}}+{{b\over e}}\to 1$$
using properties 2-3).  Property 1$'$) follows from this and property 1).  Finally, since
$${{bf}\over {xc}}={{(ac+bd)b}\over{xce}}={b\over e}{a\over x} + {{b^2}\over{e^2}}{{de}\over{xc}}\to\infty$$
using the Ptolemy equation and properties 2-3), property 2$'$) holds as well.

Applying these two types of Whitehead moves along edges in $K$, we may alter $G'$ to arrange that the edge-path for $K$ in the resulting graph
makes exactly one left turn.  Furthermore as we have just proved, properties 1-3) continue to hold for the resulting graph.

We shall complete the proof by calculating that the absolute value of the trace of the holonomy of the edge-path $K$
is asymptotic to 2, and the second key point about properties 1-3) is that they are 
sufficient to guarantee this.  To this end, let us adopt the notation that $K$ traverses the consecutive
edges $y_1,\ldots ,y_{n+1}$, the unique half-edge on the right is $x_0$, which is incident on the common
endpoint of $y_{n+1},y_1$, and the consecutive half-edges on the left are $x_1,\ldots ,x_n$, where $x_k$
has common endpoint in $K$ with $y_k,y_{k+1}$, for $k=1,\ldots,n$.  As usual identifying an edge 
or a half-edge with its lambda length, which depends upon the parameter $t$, let us define
$$
\zeta_1^2={{y_2y_{n+1}}\over{x_1x_0}}, \zeta_{n+1}^2={{x_0x_n}\over{y_1y_n}},
~{\rm and}~\zeta_k^2={{y_{k+1}x_{k-1}}\over{x_k y_{k-1}}},~{\rm for}~k=2,\ldots ,n,
$$
so the cross ratio of edge $y_k$ is $\zeta_k^{-2}$, for $k=1,\ldots ,n+1$.
The path-ordered product of matrices to compute the holonomy of $K$ beginning from the unique
left turn is given by
$$\aligned
L&\left (\begin{array} {cc} 0&\zeta_1\\ -\zeta_1^{-1}&0\\\end{array}\right )
R\left (\begin{array} {cc} 0&\zeta_2\\ -\zeta_2^{-1}&0\\\end{array}\right )\cdots
R\left (\begin{array} {cc} 0&\zeta_{n+1}\\ -\zeta_{n+1}^{-1}&0\\\end{array}\right )\\
\\
&=\left (\begin{array} {cc} ~~\zeta_1^{-1}&0\\ -\zeta_1^{-1}&\zeta_1\\\end{array}\right )
\left (\begin{array} {cc} \zeta_2^{-1}&-\zeta_2\\ 0&~~\zeta_2\\\end{array}\right )\cdots
\left (\begin{array} {cc} \zeta_{n+1}^{-1}&-\zeta_{n+1}\\ 0&~~\zeta_{n+1}\\\end{array}\right )\\
\\
&=\left (\begin{array} {cc} ~~\zeta_1^{-1}&0\\ -\zeta_1^{-1}&\zeta_1\\\end{array}\right )
\left (\begin{array} {cc} ~~(\zeta_2\cdots\zeta_{n+1})^{-1}&
\zeta_2\cdots\zeta_{n+1}\sum_{k=1}^n\prod_{j=2}^k \zeta_j^{-2}\\ 0&\zeta_2\cdots\zeta_{n+1}\\\end{array}\right ),\\
\endaligned$$
so the trace is found to be
$$(\zeta_1\zeta_2\cdots \zeta_{n+1})+ (\zeta_1\zeta_2\cdots \zeta_{n+1})^{-1}
-\zeta_1^{-2}(\zeta_1\zeta_2\cdots \zeta_{n+1})\sum_{k=1}^n\prod_{j=2}^k \zeta_j^{-2}.$$
Finally, direct calculation shows that the product telescopes, and
$$(\zeta_1\zeta_2\cdots \zeta_{n+1})={{y_{n+1}}\over{y_1}}\to 1$$
since ${{y_1}\over {y_{n+1}}}\to 1$ by property 3).  Furthermore, 
$$\zeta_1^{-2}={{x_0x_1}\over{y_2y_{n+1}}}\sim {{x_0x_1}\over{y_1y_2}}\to 0$$
by properties 2-3), and indeed, the general term in the sum also telescopes
$$\zeta_1^{-2}\zeta_2^{-2}\cdots \zeta_k^{-2}\sim {{x_0x_k}\over{y_ky_{k+1}}}\to 0,~{\rm for}~
k=2,\ldots ,n,$$
again by properties 2-3).  The absolute value of the trace is thus indeed asymptotic to 2.
Since the absolute value of the trace is twice the hyperbolic cosine of half the hyperbolic length,
the curve $K$ is asymptotically
short as $t\to\infty$.

For the final assertion, we must show that the edge-path of an essential short curve $K$ for $\lambda_t$  lies in $\partial {\mathcal A}$.  We shall use the Collar Lemma \cite{Buser}
that an essential simple closed curve of hyperbolic length $\ell$ has an embedded collar of width at least the logarithm of the hyperbolic
cotangent of ${\ell\over 2}$.  Since the dual of $G$ is an ideal triangulation and $K$ is essential, it thus follows that if $K$ traverses an edge $e$ of $G$, then the lambda length of (the ideal arc dual to) $e$ has 
divergent lambda length.  It follows that $K$ shares an edge with $E^k$, for some $k\geq 1$.
Since two essential short curves cannot intersect, again by the Collar Lemma, we conclude that
$K$ cannot meet $\partial{\mathcal A}(\lambda _t)$ for large $t$, so
in fact the edge-path for $K$ is contained in $E^k-E^{k+1}$.  

If $K$ is not homotopic to a boundary component of $G_{E^k}$, then its edge-path must make both right and left turns in $G_{E^k}$.  Without loss, we may assume that there is a left turn followed by a right turn and adopt the following notation.
Suppose that the edge-path for $K$ traverses edges $y_0,\ldots ,y_{n+1}$ in $K$
with half-edge $x_j$ incident on the common endpoint of $y_j,y_{j+1}$ for $j=1,\ldots , n$, where
$x_0\in E^k$ lies on the right and $x_n$ lies on the left of $K$.  In the extreme case that $n=1$,
the dual arcs to $x_0,y_1,x_1,y_0\in E^k$ are the consecutive edges of an ideal quadrilateral whose cross ratio is bounded near one by Lemma~\ref{horos} since the lambda lengths $x_0,x_1,y_0, y_1$ are comparable, i.e., the limit of the ratio of any pair is finite and non-zero.  The arcs dual to $y_0$ and $y_1$ are therefore a bounded distance apart, contradicting that $K$ is short.  This extreme
case gives a lower bound to the distance between the arcs dual to $y_0$ and $y_n$, so in any case,
$K$ cannot be short.  This contradiction establishes the final assertion and completes the proof.
 \end{proof}~
 
~~ \vskip .2in

\section{Proof of main result}

\begin{theorem}
The cell $C(G)$ in decorated Teichm\"uller space corresponding
to the fatgraph $G$ is asymptotic to a stable curve with pinch curves
$K$ if and only if $K$ is homotopic to the collection of edge-paths
$\partial {\mathcal A}$ for some screen ${\mathcal A}$ on $G$.
\end{theorem}

\begin{proof}
Suppose that $\lambda_t\in{\mathbb R}^0_{>0}$ is a path of lambda lengths in $C(G)$ whose projectivization
$\bar\lambda _t$ accumulates at
some point of $P({\mathbb R}^0_{\geq0})$.  Since $C(G)$ is path connected, there is a stable
path, still denoted $\lambda_t$, whose limit point is this accumulation point.  By Lemma \ref{shortcurves}, the short curves for this limit point are the components represented by edge-paths in
$\partial{\mathcal A}(\lambda _t)$.  

Conversely for any trivalent fatgraph $G$ and any screen ${\mathcal A}$ on $G$, Lemma~\ref{cubicrealization} shows that 
$\partial {\mathcal A}$ is realized as the set of short curves for a stable path in $C(G)$.

For a general not necessarily trivalent fatgraph, we require a further ingredient for the converse, namely:

\vskip .2in

\begin{theorem}\label{cyclic} ~\cite{Penner87}~
For any cyclically ordered tuple $x_1,\ldots ,x_n$ of positive real numbers satisfying the 
generalized strict triangle inequalities $x_j< \sum _{i\neq j} x_i$, for $j=1,\ldots ,n$, there is a 
cyclic Euclidean planar polygon (i.e., the polygon inscribes in a circle) unique up to orientation-preserving isometry of the plane which realizes these numbers as its consecutive edge lengths.  
\end{theorem}

\vskip .2in

To apply this result, let $L^+$ denote the collection of isotropic vectors in
Minkowski space with positive $z$-coordinate.  Given a collection of coplanar points
in $L^+$ lying in an affine plane determining an elliptic conic section, we may apply a Minkowski isometry to arrange that the plane containing these points is horizontal.  The restriction of the Minkowski pairing to this horizontal plane is a multiple of the Euclidean metric induced on the plane, so the projectivized lambda lengths of pairs of these points agree with the projectivized Euclidean lengths in the horizontal plane.  Furthermore, the intersection of the horizontal plane with $L^+$ is a round circle in this Euclidean structure.

Now, given any fatgraph $G'$ with vertices at least trivalent and any screen ${\mathcal A}'$ on $G'$, again define lambda lengths on the edges of $G'$ by $\lambda _t'(e)=t^{d_e}$, where $d_e$ is the depth of $e$ in ${\mathcal A}'$.  According to Theorem~\ref{cyclic}, the previous paragraph, 
and Theorem~\ref{lambda}, this does indeed determine a path in ${C}(G')$.

Choose any trivalent fatgraph $G$ which collapses to $G'$.  The lambda lengths on the edges
of $G-G'$ are thus given by the Euclidean lengths of the corresponding cyclic polygon again according to Theorem~\ref{cyclic}, the lambda lengths on the edges of $G'$ have already been specified, so $\lambda _t$ is now determined on $G$.

It follows that there is a unique screen ${\mathcal A}$ on $G$
which restricts to ${\mathcal A}'$ in the natural sense with corresponding lambda lengths
$\lambda _t$ on $G$, and satisfying $\partial {\mathcal A}'=\partial {\mathcal A}$ as homotopy
classes.  The proof of Lemma \ref{shortcurves} applies to $\lambda_t$ to conclude that the components of $\partial{\mathcal A}$ are precisely the short curves for $\lambda_t$.
\end{proof}

\section{Closing remarks}

In sketch, we have seen that simplicial coordinates are given explicitly in terms of lambda lengths, 
and the no vanishing cycle condition is necessary and sufficient to guarantee that these formulae are uniquely invertible.  Writing the inverse explicitly is the basic ``arithmetic problem" in decorated Teichm\"uller theory \cite{Penner95}.

One ingredient, which is related to the asymptotics of this arithmetic problem, towards describing the Deligne-Mumford compactification is:

~~\vskip .2in

\begin{theorem}~\cite{Penner96}
~\label{IJ}  Suppose that $\lambda _t$ is a stable one-parameter family 
of lambda lengths on the fatgraph $G$ with no vanishing cycles of corresponding
simplicial coordinates $X_t$.
Define $I=\{e\in E:\lambda_t(e)\to\infty\}$ and
$J=\{e\in E:X_t(e)\to 0\}$.  Then $I\subseteq J$ and $R(G_J)=G_I$,
where $R(X)$ denotes the maximal recurrent subset of $X$.
\end{theorem}

\vskip .2in

This result in tandem with Theorem~\ref{mainthm} has interesting consequences:  Take a straight-line
path in the natural affine structure of simplicial coordinates on ${C}(G)$ for some fatgraph $G$
which limits to a point that fails to satisfy the no vanishing cycle condition.  Let $E_1\subseteq E(G)$
denote the subset of edges of $G$ whose simplicial coordinates vanish, and let $R_1\subseteq E_1$ denote its maximal recurrent subset.  Depending upon the affine path, certain lambda lengths of edges in $R_1$ diverge at various rates, i.e., a screen magically pops out as determined by the arithmetic problem.

The set of all screens on all fatgraphs of a fixed topological type $g,s$ forms a partially ordered set in the natural way, where the face relation is generated as usual by collapsing edges and now also by
inclusion of screens.    The mapping class group $MC(F_g^s)$ acts on the geometric realization of this screen poset in the natural way.  As will be described in a forthcoming paper by the first-named author, a quotient of this
screen poset  is homotopy equivalent to the Deligne-Mumford compactification $\bar{\mathcal M}(F_g^s)$, where the equivalence relation is finer than simply taking $MC(F_g^s)$-orbits; more specifically,
the co-operad structure of $\bar{\mathcal M}(F_g^s)$ can be described explicitly using ``partially decorated surfaces'' (cf. \cite{Penner87, Penner06}) leading to the usual graphical diagrams of irreducible components for stable curves, and the finer equivalence relation is generated also by combinatorial isomorphisms of these diagrams.

\bibliographystyle{amsplain}

\end{document}